\documentclass[11pt]{article} \makeatother
\usepackage{amsfonts,amsmath,amssymb}
\allowdisplaybreaks[3]

\begin{document}

\numberwithin{equation}{section}
\newtheorem{theorem}{\ \ \ \ Theorem}[section]
\newtheorem{proposition}[theorem]{\ \ \ \ Proposition}
\newtheorem{lemma}[theorem]{\ \ \ \ Lemma}
\newtheorem{remark}{\ \ \ \ Remark}
\newcommand{\be}{\begin{equation}}
\newcommand{\ee}{\end{equation}}
\newcommand\bes{\begin{eqnarray}} \newcommand\ees{\end{eqnarray}}
\newcommand{\bess}{\begin{eqnarray*}}
\newcommand{\eess}{\end{eqnarray*}}
\newcommand\ds{\displaystyle}

 \begin{center} {\bf\Large On decay and blow-up of solutions for  a singular  nonlocal  }\\[2mm]
  {\bf\Large viscoelastic problem  with a nonlinear source term}
\\[4mm]
  {\large Wenjun Liu, Yun Sun  and Gang Li}\\[1mm]
{\small  College of Mathematics and Statistics, Nanjing University of
Information Science and Technology, Nanjing 210044, China.} \\[2mm]
\end{center}

\setlength{\baselineskip}{17pt}{\setlength\arraycolsep{2pt}

\begin{quote}
\noindent {\bf Abstract:} In this paper we consider a singular nonlocal viscoelastic problem with a nonlinear source term and a possible damping term. We proved that if the initial data enter into the stable set,  the solution exists globally and decays to zero with a more general rate, and if the initial data
enter into the unstable set, the solution with non-positive  initial
energy as well as  positive initial energy blows up in finite time.
These are achieved by using  the potential well theory, the modified convexity method and the perturbed
energy method.


\end{quote}

\setlength{\baselineskip}{17pt}{\setlength\arraycolsep{2pt}

\section{Introduction}

In this paper,  we investigate the following one-dimensional viscoelastic  problem with a nonlocal boundary condition
 \bes\left\{\begin{array}{ll}
\displaystyle  u_{tt}-\frac{1}{x}(x u_{x})_{x}+\int_{0}^{t}g(t-s)\frac{1}{x}(x u_{x}(x, s))_{x}{\rm d}s+au_{t}=|u|^{p-2}u,  &x\in (0, \ell), t\in (0, \infty), \medskip\\\medskip
  \displaystyle u(\ell, t)=0, \int_{0}^{\ell}x u(x, t){\rm d}x=0 &t\in [0, \infty), \\\medskip
\displaystyle u(x, 0)=u_{0}(x), \quad u_{t}(x, 0)=u_{1}(x), \quad
&x\in [0, \ell],
 \end{array}\right.\label{1.1}
 \ees
where $a\geq0$, $\ell<\infty$, $p>2$ and
$g:\mathbb{R^{+}}\rightarrow \mathbb{R^{+}}$.

This type of evolution problems, with nonlocal constraints, are generally encountered in heart transmission theory, thermoelasticity,
chemical engineering, underground water flow, and plasma physics. The nonlocal boundary conditions arise mainly when the data on the boundary can not be measured directly, but their average values are known.
We can refer to the works of Cahlon and Shi \cite{ckp95}, Cannon \cite{c1963}, Choi and Chan  \cite{ck1992}, Ewing and Lin \cite{el1991}, Ionkin \cite{i1977}, Kamynin \cite{k1964}, Samarskii \cite{s1980}, and Shi and Shilor \cite{ss1990}. The first paper discussed second order partial differential equations with nonlocal integral conditions goes back to Cannon \cite{c1963}.
In fact, most of the works were about the classical solutions. Later, Mixed problems with classical and nonlocal (integral) boundary conditions related to parabolic and hyperbolic equations have been extensively established and several results concerning existence and uniqueness have been considered by  Bouziani \cite{b1996}, Ionkin \cite{i1977}, Kamynin \cite{k1964}, Mesloub \cite{mm2002}, Pulkina \cite{p1999}.

In the absence of the viscoelastic term (i.e., $g=0$), Mesloub  and
Bouziani \cite{mb1999} studied the following equation
$$v_{tt}-\frac{1}{x}v_{x}-v_{xx}=f(x, t),\quad x\in(0, \ell), \quad t\in(0, T), $$
and obtained the  existence and uniqueness of a strong solution. Later, Mesloub and Messaoudi  \cite{mm2002} solved a three-point boundary-value problem
for a hyperbolic equation with a Bessel operator and an integral condition based on an energy method. Then in \cite{mm2003} they considered a nonlinear
one-dimensional hyperbolic problem with a linear damping term and established a blow-up result for large initial data and a decay result for small initial data.

In the presence of the viscoelastic term (i.e., $g\neq0$), Mecheri et
al. \cite{mss9} studied the following equation
$$u_{tt}-\frac{1}{x}(x u_{x})_{x}+\int_{0}^{t}g(t-s)\frac{1}{x}(x u_{x}(x, s))_{x}{\rm d}s+au_{t}=f(x,t),\quad 0<x<1,\quad t>0,$$
for $a>0$ and proved the existence and uniqueness of the strong solution.  Then, Mesloub et al. \cite{mm2010} considered a nonlinear mixed problem for a viscoelastic equation with a dissipation term under a nonlocal boundary condition and obtained the existence and uniqueness of the weak solution based on the iteration processes. Later, the global existence, decay and blow-up of solutions of problem \eqref{1.1}
(when $a=0$) were established by Mesloub and Messaoudi in \cite{mm10}, where the authors studied the blow-up result with only negative initial energy. Recently, Wu \cite{wu2011} improved \cite{mm10} by establishing the blow-up result  with nonpositive initial energy as well as positive initial energy.


For the case of initial and boundary value problems for linear and nonlinear viscoelastic equations with classical conditions, many results have also been extensively studied. Cavalcanti et al. \cite{cvs2002} studied}
$$u_{tt}-\Delta u+\int_{0}^{t}g(t-\tau)\Delta u(\tau){\rm d}\tau+a(x)u_{t}+|u|^m u=0, \quad (x, t)\in \Omega\times(0, \infty), $$
for $a:\Omega\rightarrow \mathbb{R^+}$, a function, which may be null on a part of the domain $\Omega$. Under the conditions that $a(x)\geq a_{0}>0$
on $\omega\subset\Omega$, with $\omega$ satisfying some geometry restrictions and $$-\xi_{1} g(t)\leq g'(t)\leq -\xi_{2} g(t), \quad t\geq0, $$
the authors established an exponential rate of decay. Berrimi and Messaoudi \cite{bm2004} improved Cavalcanti's result by introducing a different functional which allowed to weak the conditions on both $a$ and $g$. In particular, the function $a(x)$ can vanish on the whole domain $\Omega$ and consequently the geometry condition has disappeared. In \cite{co2003}, Cavalcanti et al. considered
$$u_{tt}-k_{0}\Delta u+\int_{0}^{t}{\rm div}[a(x)g(t-\tau)\nabla u(\tau)]{\rm d}\tau+b(x)h(u_{t})+f(u)=0, $$
under similar conditions on the relaxation function $g$ and
$a(x)+b(x)\geq\rho>0$, for all $x\in\Omega$. They improved the
result of \cite{cvs2002} by establishing exponential stability for
$g$ decaying exponentially and $h$ linear and polynomial stability
for $g$ decaying polynomially and $h$ nonlinear. In \cite{bm06},
Berrimi and  Messaoudi considered
$$u_{tt}-\Delta u+\int_{0}^{t}g(t-\tau)\Delta u(\tau){\rm d}\tau=|u|^{p-2}u$$
in a bounded domain and $p>2$. They established a local existence
result and showed that, under weaker condition $g'(t)\leq\xi
g^{r}(t)$, the solution is global and decay in a polynomial or
exponential fashion when the initial data is small enough. Then
Messaoudi \cite{m2008} improved this result by establishing a
general decay of energy which is similar to the relaxation function
under weaker condition that $g'(t)\leq\xi(t) g(t)$. In regard of
nonexistence, Messaoudi \cite{am2003} considered
$$u_{tt}-\Delta u+\int_{0}^{t}g(t-\tau)\Delta u(\tau){\rm d}\tau+a|u_{t}|^{m-2}u_{t}=|u|^{p-2}u$$
and established a blow up result for solutions with negative energy if $p>m$ and a global existence result for $p\leq m$. Then Messaoudi \cite{am2006}
improved this result by accommodating certain solutions with positive initial energy. Liu \cite{l10} obtained the similar blow-up result  for
the viscoelastic problem with strong damping and nonlinear source by using  the potential
well theory and convexity technique. For other related works, we refer the readers to \cite{l73, ly11,mn03,am9,sw2006,ts2012,y2009,z2011,yz2005} and the references therein.

Inspired by \cite{bm06, l10, mm10, m2008}, we intend to study the
blow-up and decay properties of problem \eqref{1.1} in this paper.
Our goal is to establish a decay result with a more general rate and a blow-up result with non-positive  initial
energy as well as  positive initial energy. The main difficulties we
encounter here arise from the simultaneous appearance of the
singular nonlocal viscoelastic term, the possible damping term, as
well as the nonlinear source term. We first show that if the initial data
enter into the unstable set, the source term is enough to obtain
blow-up result no matter $a=0$ or $a>0$. This is achieved by using the potential well theory and the modified convexity method.
We then establish the decay
result under the condition that $g'(t)\leq - \xi(t)g^{r}(t)$, which
is more general than that of \cite{bm06, m2008}, by constructing
some functionals and using the perturbed
energy method.

The paper is organized as follows. In Section \ref{s2} we present some assumptions and known results and state the main results.
  Section \ref{s3} is devoted to the proof of  the blow-up  result. The decay result is proved in Sections \ref{s4}.

\section{Preliminaries and main results}\label{s2}

In this section we first introduce some functional spaces and present some assumptions and known results which will be used throughout this work.

Let $L^{p}_{x}=L^{p}_{x}(0, \ell)$ be the weighted Banach space equipped with the norm $$\|u\|_{p}=\left(\int_{0}^{\ell}x|u|^{p}{\rm d}x\right)^{\frac{1}{p}}. $$
In particular, when $p=2$, we denote $H=L^{2}_{x}(0,\ell)$ to be the weighted Hilbert space of square integrable functions having the finite norm$$\|u\|_{H}=\left(\int_{0}^{\ell}xu^{2}{\rm d}x\right)^{\frac{1}{2}}. $$
We take $V=V^{1, 1}_{x}(0, \ell)$ to be the weighted Hilbert space equipped with the norm$$\|u\|_{V}=\left(\|u\|_{H}^{2}+\|u_{x}\|_{H}^{2}\right)^{\frac{1}{2}}, $$
and $$V_{0}=\{u\in V \ \text{such that}\ u(\ell)=0\}.$$

For the relaxation function $g$, we give the following assumptions:
 \begin{quote}
(G1) $g(t):\mathbb{R^{+}}\rightarrow \mathbb{R^{+}}$ is a
non-increasing $C^2$ function such that
\begin{equation*}
 g(0)>0, \qquad 1-\int_{0}^{\infty }g(s){\rm d}s = l>0.
\end{equation*}
(G2) There exists a positive differentiable function $\xi(t)$ such
that
\begin{equation}\label{2.30}
  g'(t)\leq - \xi(t)g^{r}(t), \quad t\geq 0, 1\leq r<\frac{3}{2},
\end{equation}
and $\xi(t)$ satisfies, for some positive constant $L$,
 $$\left|\frac{\xi '(t)}{\xi(t)}\right|\leq L, \quad \xi'(t)\leq 0,\quad \int_{0}^{+\infty}\xi(s){\rm d}s=+\infty, \quad \forall\ t>0.$$
Furthermore, when $1<r<\frac{3}{2}$, for any fixed $t_{0}>0$, there
exists a positive constant $C_{r}$ depending only on $r$, such that
\begin{equation}\label{1.1.1.2}
\frac{t}{\left(1+\int_{t_{0}}^{t}\xi(s){\rm
d}s\right)^\frac{1}{2(r-1)}}\leq C_{r},\quad \forall\ t\geq
t_{0}.
\end{equation}
\end{quote}

\begin{remark} \label{re1} Condition $r<\frac{3}{2}$ is made to ensure that $\int_{0}^{\infty}g^{2-r}(s){\rm d}s<\infty.$
\end{remark}
\begin{remark} \label{re2}
If $\xi(t)\equiv\xi=$contant, {\rm (G2)} recaptures that of \cite{bm06,l10, mm10}. If $r\equiv1$, {\rm (G2)} recaptures  that of \cite{m2008,am9}. Therefore, {\rm (G2)} is a generalization of \cite{bm06,l10, mm10,m2008,am9}. In  particular, when  $\xi(t)\equiv\xi$ and $1< r<\frac{3}{2}$, \eqref{1.1.1.2} holds naturally. 
\end{remark}

\begin{lemma}   \label{le2.1} {\rm(\cite{mm10}, Poincar\'{e}-type inequality)}
For any $v$ in $V_{0}$, we have$$\int_{0}^{\ell}x v^{2}(x){\rm d}x\leq C_{p}\int_{0}^{\ell}xv_{x}^{2}(x){\rm d}x,$$
where $C_{p}$ is some positive constant.
\end{lemma}

\begin{lemma}   \label{le2.2}{\rm(\cite{mm10})}
For any $v$ in $V_{0}$, $2<p<4$, we have$$\int_{0}^{\ell}x |v|^{p}{\rm d}x\leq C_{*}\|v_{x}\|_{2}^{p}, $$ where $C_{*}$ is a constant depending on $\ell$ and $p$  only.
\end{lemma}

We state, without proof, a local existence result for problem \eqref{1.1}.
The proof can be easily established by adopting the arguments of \cite{bm06}, \cite{mm2010} and  \cite{mm2003}.

\begin{theorem}\label{th2.3}
Suppose that {\rm (G1)} holds and  $2<p<3$. Then for any $u_{0}$ in
$V_{0}$ and $u_{1}$ in $H$, problem \eqref{1.1} has a unique local
solution
$$u\in  C(0, T_{max};V_{0})\cap C^{1}(0,T_{max};H)$$ for $T_{max}>0$ small enough.
\end{theorem}

\begin{remark} \label{re4}
The condition $2<p<3$ is needed so that the embedding of $V_{0}$ in $L_{x}^{2}$ is Lipschit {\rm(see \cite[Lemma 5.2]{mm2003})}.
\end{remark}

Next we introduce the functionals for $I(t)$, $J(t)$ and $E(t)$:
\bes
I(t)&:=&I(u(t))=\Big(1-\int_{0}^{t}g(s)\, {\rm d}s\Big)\int_{0}^{\ell}x u_{x}^{2}{\rm d}x  + (g\circ u_{x})(t)-\int_{0}^{\ell}x|u(t)|^p{\rm
d}x, \label{2.1}\\
J(t)&:=&J(u(t)) =\frac{1}{2}\Big(1-\int_{0}^{t}g(s)\, {\rm d}s\Big)\int_{0}^{\ell}x u_{x}^{2}{\rm d}x  +\frac{1}{2} (g\circ u_{x})(t)-\frac{1}{p}\int_{0}^{\ell}x|u(t)|^p{\rm
d}x, \label{2.2}\\
E(t)&:=&E(u(t))=J(t)+\frac{1}{2}\int_{0}^{\ell}x u_{t}^{2}{\rm d}x,  \label{2.3}
 \ees
where $$(g\circ u_{x})(t)=\int_{0}^{\ell}\int_{0}^{t}xg(t-s)|u_{x}(x, t)-u_{x}(x, s)|^{2}{\rm d}s{\rm d}x.$$

\begin{remark} \label{re3}
 A multiplication of equation \eqref{1.1} by $xu_{t}$ and integration over $(0, \ell)$ easily yields
\begin{equation}\label{2.4}
E'(t)=\frac{1}{2}(g'\circ u_{x})(t)-\frac{1}{2}g(t)\int_{0}^{\ell}x u_{x}^{2}{\rm d}x-a\int_{0}^{\ell}x u_{t}^{2}{\rm d}x\leq-a\int_{0}^{\ell}x u_{t}^{2}{\rm d}x\leq0, \quad \forall\ t\geq0.
\end{equation}
\end{remark}

We are now in a position to state our main results.

\begin{theorem}\label{th2.5}
Assume that {\rm (G1)} holds and $2<p<3$, let $u$ be the unique
local solution to problem \eqref{1.1} and denote
$d_1=\frac{p-2}{2p}\left(\frac{l}{C_{*}^{2/p}}\right)^{\frac{p}{p-2}}$.
For any fixed $\delta<1$, assume that $u_{0},$ $u_{1}$ satisfy
\begin{equation}\label{2.6}
E(0)<\delta d_{1}, \qquad I(0)<0.
\end{equation}
Suppose that
\begin{equation}\label{2.7}
 \int_0^{\infty} g(s){\rm d}s\leq\frac{p-2}{p-2+1/[(1-\hat{\delta})^2p+2\delta(1-\hat{\delta})]}
\end{equation}
where $\hat{\delta}=\max\{0, \delta\}$. Then the solution of problem
\eqref{1.1} blows up in a finite time $T^*$ in the sense that
\begin{equation*}
\lim_{t\rightarrow T^{*-}}{\|u\|_{H}^2}=+\infty.
\end{equation*}.
\end{theorem}

\begin{theorem}\label{th2.4}
Assume that {\rm (G1)} holds and $2<p<3$, let $u$ be the unique
local solution to problem \eqref{1.1}. In addition, assume that
$u_{0},$ $u_{1}$ satisfy
\begin{equation}
E(0)<d_1,\qquad I(0)>0.
\end{equation}\label{2.5}
Then the solution $u$ is global and satisfies
\begin{equation}\label{2.5'}
\int_{0}^{\ell}x u_{x}^{2}{\rm
d}x\leq\frac{2p}{l(p-2)}E(t)\leq\frac{2p}{l(p-2)}E(0),\quad \forall\
t>0.
\end{equation}
\end{theorem}

\begin{theorem}\label{th2.6}
Under the assumptions of theorem \ref{th2.4}, suppose further that {\rm (G2)} holds. Then for each $t_{0}>0$, there exist positive constants $K$ and $\kappa$ such that
\begin{equation}
E(t)\leq
\left\{
 \begin{aligned}
&Ke^{-\kappa\int_{t_0}^t\xi(s){\rm d}s}, &&r=1, \\
&K\left(1+\int_{t_0}^{t}\xi(s){\rm d}s\right)^{-\frac{1}{r-1}}, &&1<r<\frac{3}{2}.
  \end{aligned}
  \right.\label{2.8}
  \end{equation}
\end{theorem}

\begin{remark}\label{re5}

Note that when $1< r<\frac{3}{2}$, we obtain more general type of
decays. If we choose $\xi(t)\equiv\xi$, \eqref{2.8} gives the
polynomial rate decay as $E(t)\leq K(1+t)^{-\frac{1}{r-1}}$, which
coincides with the results of \cite{bm06,l10, mm10}.  If we choose
$\xi(t)=(1+t)^{-m}$ for $0<m<3-2r<1$ $($ which satisfies
\eqref{1.1.1.2}$)$, we have $g(t)\leq \frac{C_{0}}{(1+t)^q}$ with
$q=\frac{1-m}{r-1}$ and \eqref{2.8}  also gives the polynomial rate
of decay as $E(t)\leq\frac{C_{1}}{(1+t)^q}$. In particular, if we
choose $\xi(t)=\frac{2(r-1)}{(t+1)^{3-2r}}+\frac{1}{t+1}$, which
satisfies {\rm (G2)}, then we have
$g(t)\leq\frac{C}{\left[(t+1)^{2(r-1)}+\ln(t+1)-1\right]^{\frac{1}{r-1}}}$
and  a new type of decay as $E(t)\leq\frac{K}{\left[(t+1)^{2(r-1)}+\ln(t+1)-1\right]^{\frac{1}{r-1}}}$ is established.
\end{remark}
\section{Blow-up of solutions}\label{s3}

In this section, we prove a finite time blow-up result for initial data in the unstable set.

For $t\ge 0$, we define
\[d(t)=\inf_{u\in  V_{0}\setminus\{0\}}\sup_{\lambda\geq 0}J(\lambda u)\]
and
\begin{equation}\label{3.09181}
\mathcal{N}=\{u\in V_{0}\setminus\{0\}:\quad I(u(t))=0\}.
\end{equation}
Then we can prove the following lemma.

\begin{lemma}\label{le3.1}
For $t\ge 0$,  we have
$$ 0<d_1\leq d(t)\leq d_2(u)=\sup_{\lambda\geq 0} J(\lambda u)$$
and\begin{equation}\label{3.09183}
d(t)=\inf_{u\in\mathcal{N}}{J(u)}.
\end{equation}
\end{lemma}

 {\bf Proof.}\ Obviously, $$d(t)\leq d_2(u)=\sup_{\lambda\geq 0} J(\lambda u).$$
 Since
 $$ J(\lambda u)=\frac{\lambda^{2}}{2}\left[\left(1-\int_{0}^{t}g(s){\rm d}s\right)\int_{0}^{\ell}x u_{x}^{2}{\rm d}x + (g\circ u_{x})(t)\right]-\frac{\lambda^{p}}{p}\int_{0}^{\ell}x|u|^{p}{\rm d} x.$$
 We get
 $$\frac{{\rm d}}{{\rm d}\lambda}J(\lambda u)=\lambda\left[\left(1-\int_{0}^{t}g(s){\rm d}s\right)\int_{0}^{\ell}x u_{x}^{2}{\rm d}x + (g\circ u_{x})(t)\right]-\lambda^{p-1}\int_{0}^{\ell}x|u|^{p}{\rm d} x.$$
 Let $$\frac{{\rm d}}{{\rm d}\lambda}J(\lambda u)=0, $$ which implies
 $$\bar{\lambda_{1}}=0, \qquad \bar{\lambda_{2}}=\left[\frac{\left(1-\int_{0}^{t}g(s){\rm d}s\right)\int_{0}^{\ell}x u_{x}^{2}{\rm d}x + (g\circ u_{x})(t)}{\int_{0}^{\ell}x|u|^{p}{\rm d} x}\right]^{\frac{1}{p-2}}.$$
 An elementary calculation shows
 $$\frac{{\rm d^2}}{{\rm d}\lambda^{2}}J(\bar{\lambda_{1}} u)>0 \qquad \text{and}\ \qquad \frac{{\rm d^2}}{{\rm d}\lambda^{2}}J(\bar{\lambda_{2}} u)<0.$$
 Using  {\rm (G1)} and Lemma \ref{le2.2}, we get
 \begin{equation*}
 \begin{aligned}
 \sup_{\lambda\geq 0} J(\lambda u)=J(\bar{\lambda_{2}} u)&=\frac{p-2}{2p}\left[\frac{\left(1-\int_{0}^{t}g(s){\rm d}s\right)\int_{0}^{\ell}x u_{x}^{2}{\rm d}x + (g\circ u_{x})(t)}{\left(\int_{0}^{\ell}x|u|^{p}{\rm d} x\right)^{2/p}}\right]^{\frac{p}{p-2}}\\
 &\geq \frac{p-2}{2p}\left[\frac{l\int_{0}^{\ell}x u_{x}^{2}{\rm d}x}{\left(\int_{0}^{\ell}x|u|^{p}{\rm d} x\right)^{2/p}}\right]^{\frac{p}{p-2}} \geq \frac{p-2}{2p}\left(\frac{l}{C_{*}^{2/p}}\right)^{\frac{p}{p-2}}\\&=d_{1}>0,
  \end{aligned}
\end{equation*}
which implies that $d(t)\geq d_{1}.$

To  get \eqref{3.09183},
straightforward computations lead to
\begin{align}\label{3.09184}
&I(\bar{\lambda_{2}}u)\nonumber\\
 =&\left(1-\int_{0}^{t}g(s){\rm
d}s\right)\int_{0}^{\ell}x\left(\bar{\lambda_{2}}u\right)_{x}^{2}{\rm
d}x+(g\circ
\left(\bar{\lambda_{2}}u\right)_{x})(t)-\int_{0}^{\ell}x|\bar{\lambda_{2}}u|^{p}{\rm
d}x\nonumber\\
=&\left[\frac{\left(1-\int_{0}^{t}g(s){\rm
d}s\right)\int_{0}^{\ell}x u_{x}^{2}{\rm d}x + (g\circ
u_{x})(t)}{\int_{0}^{\ell}x|u|^{p}{\rm d}
x}\right]^{\frac{2}{p-2}}\left[\left(1-\int_{0}^{t}g(s){\rm
d}s\right)\int_{0}^{\ell}x u_{x}^{2}{\rm d}x+
(g\circ u_{x})(t)\right]\nonumber\\
&-\left[\frac{\left(1-\int_{0}^{t}g(s){\rm
d}s\right)\int_{0}^{\ell}x u_{x}^{2}{\rm d}x + (g\circ
u_{x})(t)}{\int_{0}^{\ell}x|u|^{p}{\rm d}
x}\right]^{\frac{p}{p-2}}\int_{0}^{\ell}x|u|^{p}{\rm d}x\nonumber\\
=&\frac{\left[\left(1-\int_{0}^{t}g(s){\rm
d}s\right)\int_{0}^{\ell}x u_{x}^{2}{\rm d}x + (g\circ
u_{x})(t)\right]^{\frac{p}{p-2}}}{\left(\int_{0}^{\ell}x|u|^{p}{\rm
d}x\right)^{\frac{2}{p-2}}}\nonumber\\
 &\times\left\{\frac{\left(1-\int_{0}^{t}g(s){\rm
d}s\right)\int_{0}^{\ell}x u_{x}^{2}{\rm d}x+ (g\circ
u_{x})(t)}{ \left(1-\int_{0}^{t}g(s){\rm
d}s\right)\int_{0}^{\ell}x u_{x}^{2}{\rm d}x + (g\circ
u_{x})(t) }-1\right\}\nonumber\\
=&0,
  \end{align}
  which implies that $\bar{\lambda_{2}}u\in \mathcal {N}$. Also, for
any $u\in \mathcal {N}$, we note that
$$\bar{\lambda_{2}}(u)=\left[\frac{\left(1-\int_{0}^{t}g(s){\rm d}s\right)\int_{0}^{\ell}x u_{x}^{2}{\rm d}x + (g\circ u_{x})(t)}{\int_{0}^{\ell}x|u|^{p}{\rm d} x}\right]^{\frac{1}{p-2}}=1.$$
Therefore we have $\bar{\lambda_{2}}(u)u=u$ for all $u\in \mathcal
{N}$. Thus we complete the proof.

\begin{lemma}\label{le3.2}
Under the same assumptions as in Theorem \ref{th2.5}, one has $I(u(t))<0$ and
\begin{equation}\label{3.1}
d_{1}<\frac{p-2}{2p}\left[\left(1-\int_{0}^{t}g(s){\rm d}s\right)\int_{0}^{\ell}x u_{x}^{2}{\rm d}x + (g\circ u_{x})(t)\right]<\frac{p-2}{2p}\int_{0}^{\ell}x|u|^{p}{\rm d}x,
\end{equation}
for all $t\in[0, T_{max}).$
\end{lemma}

{\bf Proof.}\ Using \eqref{2.4} and {\rm\eqref{2.6}}, we have $E(t)\leq\delta d_{1}$ for all $t\in[0,T_{max})$. Furthermore, we can obtain $I(u(t))<0$ for all $t\in[0,T_{max})$.

In fact, if it is not true, then there exists some $t_{0}\in[0,
T_{max})$ such that $I(t_{0})\geq0$. Since  $I(0)<0$, it follows
that there exists some $\tilde{t}\in(0,t_{0}]$ such that
$I(u(\tilde{t})=0$. Define
\begin{equation}\label{3.09171}
t^*=\inf\left\{\tilde{t}\in(0,t_{0}]:\quad
\left(1-\int_{0}^{\tilde{t}}g(s){\rm d}s\right)\int_{0}^{\ell}x
u_{x}^{2}(\tilde{t}){\rm d}x+(g\circ
u_{x})(\tilde{t})=\int_{0}^{\ell}x|u(\tilde{t})|^{p}{\rm
d}x\right\}.
\end{equation}
Then, we have $I(u(t^*)=0$ and 
\begin{equation}\label{3.2}
\left(1-\int_{0}^{t}g(s){\rm d}s\right)\int_{0}^{\ell}x u_{x}^{2}{\rm d}x + \left(g\circ u_{x}\right)(t)<\int_{0}^{\ell}x|u|^{p}{\rm d}x, \quad 0\leq t<t^{*}.
\end{equation}
Next, we consider two cases:

 Case 1 : Suppose that  $\| u(t^*)\|_{H}^{2}=0$, using the
regularity of $u(t)$, we have \begin{equation}\label{3.09031}
\lim_{t\rightarrow t^{*-}}{\| u(t)\|_{H}^{2}}=0.
\end{equation}
On the other hand, from \eqref{3.2} and Lemma \ref{le2.2}, we obtain
\begin{equation}\label{3.2'}
\left(1-\int_{0}^{t}g(s){\rm d}s\right)\int_{0}^{\ell}x u_{x}^{2}{\rm d}x + \left(g\circ u_{x}\right)(t)<\int_{0}^{\ell}x|u|^{p}{\rm d}x\leq C_{*}\|u_{x}\|_{2}^p, \quad 0\leq t<t^{*},
\end{equation}
 and $\| u(t)\|_{H}^{2}\neq0$, for all $t\in[0,t^*)$. Therefore we have
  \begin{equation*} \lim_{t\rightarrow t^{*-}}{\|
u(t)\|_{H}^{2}}>\left(\frac{l}{C_{*}}\right)^{\frac{1}{p-2}},
\end{equation*}
which contradicts to \eqref{3.09031}.

Case 2 : Suppose that $\| u(t^*)\|_{H}^{2}\neq0$. Applying Lemma
\ref{le3.1}, we see that $d(t)$ is the infimum of $J(u(t))$ over all
functions $u$ in $\mathcal {N}$ and $J(u(t^*))\geq d(t)\geq d_{1}$,
which contradicts to $J(u(t^*))\leq E(t^*)<d_{1}$. Thus, we conclude
that $I(t)<0$ for all $t\in[0,T_{max})$.

To get \eqref{3.1},  we use \eqref{3.2}, Lemma \ref{le3.1} and the
conclusion that  $I(t)<0$ for all $t\in[0,T_{max})$ and get
\begin{align}
d_{1}&\leq\frac{p-2}{2p}\frac{\left[\left(1-\int_{0}^{t}g(s){\rm d}s\right)\int_{0}^{\ell}x u_{x}^{2}{\rm d}x + (g\circ u_{x})(t)\right]\left[\left(1-\int_{0}^{t}g(s){\rm d}s\right)\int_{0}^{\ell}x u_{x}^{2}{\rm d}x + (g\circ u_{x})(t)\right]^{\frac{2}{p-2}}}{\left(\int_{0}^{\ell}x|u|^{p}{\rm d}x\right)^{\frac{2}{p-2}}}\nonumber\\
&<\frac{p-2}{2p}\left[\left(1-\int_{0}^{t}g(s){\rm
d}s\right)\int_{0}^{\ell}x u_{x}^{2}{\rm d}x + (g\circ
u_{x})(t)\right], \quad 0\leq t<T_{max}.\label{3.3}
\end{align}
It follows from {\rm\eqref{3.2}} and {\rm\eqref{3.3}} that
\begin{equation*}
0<d_{1}<\frac{p-2}{2p}\left[\left(1-\int_{0}^{t}g(s){\rm
d}s\right)\int_{0}^{\ell}x u_{x}^{2}{\rm d}x + (g\circ
u_{x})(t)\right]<\frac{p-2}{2p}\int_{0}^{\ell}x|u|^{p}{\rm d}x,\quad
0\leq t<T_{max}.
\end{equation*}
Thus,  we complete the proof.

\begin{lemma}\label{le3.3}{\rm(\cite{l73})}
Let $L(t)$ be a positive $C^{2}$ function, which satisfies, for $t>0$, the inequality $$L(t)L''(t)-(1+\zeta)L'(t)^{2}\geq0$$
with some $\zeta>0$. If $L(0)>0$ and $L'(0)>0$, then there exists a time $ T^*\leq\frac{L(0)}{\zeta L'(0)}$ such that
\begin{displaymath}
\lim_{t\rightarrow T^{*-}}L(t)=\infty.
\end{displaymath}
\end{lemma}

{\bf Proof of Theorem \ref{th2.5}.} \ Assume by contradiction that the solution $u$ is global. Then, we consider $L:[0, T]\rightarrow \mathbb{R_{+}}$ defined by
\begin{equation}\label{5.14}
L(t)=\int_{0}^{\ell}xu^2{\rm d}x+a\int_{0}^{t}\int_{0}^{\ell}xu^2{\rm d}x{\rm d}s+a(T-t)\int_{0}^{\ell}xu^2_{0}{\rm d}x+b(t+T_{0})^2,
\end{equation}
where $T$, $b$ and $T_{0}$ are positive constants to be chosen later. Then $L(0)>0$. Furthermore,
\label{3.4}
\begin{align}
L'(t)&=2\int_{0}^{\ell}x u u_{t}{\rm d}x+a\int_{0}^{\ell}x(u^{2}-u_{0}^{2}){\rm d}x+2b(t+T_{0})\nonumber\\
&=2\int_{0}^{\ell}x u u_{t}{\rm d}x+2a\int_{0}^{t}\int_{0}^{\ell}x u u_{s}{\rm d}x{\rm d}s+2b(t+T_{0}),
\end{align}
and, consequently,
$$L''(t)=2\int_{0}^{\ell}x u u_{tt}{\rm d}x+2\int_{0}^{\ell}xu_{t}^{2}{\rm d}x+2a\int_{0}^{\ell}x u u_{t}{\rm d}x+2b$$
for almost every $t\in[0, T].$ Testing the equation {\rm\eqref{1.1}}} with $xu$ and plugging the result into the expression of $L''(t)$, we obtain
\begin{align*}
L''(t)=&-2\int_{0}^{\ell}x u_{x}^2{\rm d}x+2\int_{0}^{\ell}\int_{0}^{t}g(t-s)x u_{x}(x, t)u_{x}(x, s){\rm d}s{\rm d}x\nonumber\\
&-2a\int_{0}^{\ell}x u u_{t}{\rm d}x+2\int_{0}^{\ell}x|u|^{p}{\rm d}x+2\int_{0}^{\ell}x u_{t}^{2}{\rm d}x+2a\int_{0}^{\ell}x u u_{t}{\rm d}x+2b\\
=&2\bigg[\int_{0}^{\ell}x u_{t}^{2}{\rm d}x-\left(1-\int_{0}^{t}g(s){\rm d}s\right)\int_{0}^{\ell}x u_{x}^2{\rm d}x\nonumber\\
&-\int_{0}^{\ell}\int_{0}^{t}g(t-s)x u_{x}(x, t)\left(u_{x}(x, t)-u_{x}(x, s)\right){\rm d}s{\rm d}x+\int_{0}^{\ell}x|u|^{p}{\rm d}x+b\bigg]
\end{align*}
for almost every $t\in [0, T]$. Therefore, we get
\begin{equation}
\begin{aligned}
L(t)L''(t)-\frac{p+2}{4}L'(t)^2=&2L(t)\Bigg[\int_{0}^{\ell}x u_{t}^{2}{\rm d}x-\left(1-\int_{0}^{t}g(s)\,
{\rm d}s\right)\int_{0}^{\ell}x u_{x}^2{\rm d}x\nonumber\\
&-\int_{0}^{\ell}\int_{0}^{t}g(t-s)x u_{x}(x, t)\left(u_{x}(x, t)-u_{x}(x, s)\right){\rm d}s{\rm d}x\nonumber\\
&+\int_{0}^{\ell}x|u|^{p}{\rm d}x+b\Bigg]+(p+2)\bigg[\eta(t)-\left(L(t)-a(T-t)\int_{0}^{\ell}xu_{0}^2{\rm d}x\right)\nonumber\\
&\left(\int_{0}^{\ell}x u_{t}^2{\rm d}x+a\int_{0}^{t}\int_{0}^{\ell}x u_{s}^2{\rm d}x{\rm d}s+b\right)\bigg],
\end{aligned}
\end{equation}
where
\begin{equation}
\begin{aligned}
\eta(t)=&\bigg(\int_{0}^{\ell}x u^2{\rm d}x+a\int_{0}^{t}\int_{0}^{\ell}x u^2{\rm d}x{\rm d}s+b(t+T_{0})^2\bigg)\bigg(\int_{0}^{\ell}x u_{t}^2{\rm d}x+a\int_{0}^{t}\int_{0}^{\ell}x u_{s}^2{\rm d}x{\rm d}s+b\bigg) \nonumber\\
&-\bigg[\int_{0}^{\ell}x uu_{t}{\rm d}x+a\int_{0}^{t}\int_{0}^{\ell}x uu_{s}{\rm d}x{\rm d}s+b(t+T_{0})\bigg]^2.
\end{aligned}
\end{equation}
Using Schwarz's inequality, we can easily get $\eta(t)\geq0$ for every $t\in[0, T]$.
As a consequence, we reach the following differential inequality
\begin{equation}\label{3.5}
L(t)L''(t)-\frac{p+2}{4}L'(t)^2\geq L(t)\Phi(t), \quad  a.e.\quad t\in[0, T],
\end{equation}
where $\Phi:[0, T]\mapsto \mathbb{R}_{+}$ is the map defined by
\begin{equation}
\begin{aligned}
\Phi(t)=&-p\int_{0}^{\ell}x u_{t}^{2}{\rm d}x-2\left(1-\int_{0}^{t}g(s){\rm d}s\right)\int_{0}^{\ell}x u_{x}^2{\rm d}x\, -a(p+2)\int_{0}^{t}\int_{0}^{\ell}xu_{s}^2{\rm d}x{\rm d}s\nonumber\\
&-2\int_{0}^{\ell}\int_{0}^{t}g(t-s)x u_{x}(x, t)\left(u_{x}(x, t)-u_{x}(x, s)\right){\rm d}s{\rm d}x\,
+2\int_{0}^{\ell}x|u|^{p}{\rm d}x-pb\\
=&-2pE(t)+p(g\circ u_{x})(t)+(p-2)\left(1-\int_{0}^{t}g(s){\rm d}s\right)\int_{0}^{\ell}x u_{x}^2{\rm d}x-pb\nonumber\\
&-2\int_{0}^{\ell}\int_{0}^{t}g(t-s)x u_{x}(x, t)\left(u_{x}(x, t)-u_{x}(x, s)\right){\rm d}s{\rm d}x-a(p+2)\int_{0}^{t}\int_{0}^{\ell}x u_{s}^2{\rm d}x{\rm d}s.
\end{aligned}
\end{equation}
By {\rm\eqref{2.4}}, for all $t\in[0, T]$ we may also write
\begin{align}
\Phi(t)\geq&-2pE(0)+p(g\circ u_{x})(t)+(p-2)\left(1-\int_{0}^{t}g(s){\rm d}s\right)\int_{0}^{\ell}x u_{x}^2{\rm d}x-pb\nonumber\\
&-2\int_{0}^{\ell}\int_{0}^{t}g(t-s)x u_{x}(x, t)\left(u_{x}(x, t)-u_{x}(x, s)\right){\rm d}s{\rm d}x\nonumber\\
&+a(p-2)\int_{0}^{t}\int_{0}^{\ell}x u_{s}^2{\rm d}x{\rm d}s.\label{3.6}
\end{align}
By using Young's inequality, we have
\begin{align}
&2\int_{0}^{\ell}\int_{0}^{t}g(t-s)x u_{x}(x, t)\left(u_{x}(x, t)-u_{x}(x, s)\right){\rm d}s{\rm d}x\nonumber\\
\leq&\frac{1}{\varepsilon}\,
\int_{0}^{t}g(s)\int_{0}^{\ell}x u_{x}^2{\rm d}s{\rm d}x+\varepsilon(g\circ u_{x})(t),\label{3.7}
\end{align}
for any $\varepsilon>0.$
Substituting {\rm\eqref{3.7}} for the fifth term of the right hand side of {\rm\eqref{3.6}}, we obtain
\begin{align}
\Phi(t)\geq&-2pE(0)+\left[(p-2)-\left(p-2+\frac{1}{\varepsilon}\right)\int_{0}^{t}g(s){\rm d}s\right]\int_{0}^{\ell}x u_{x}^2{\rm d}x\,\nonumber
\\
&+(p-\varepsilon)(g\circ u_{x})(t)+a(p-2)\int_{0}^{t}\int_{0}^{\ell}x u_{s}^2{\rm d}x{\rm d}s-pb.\label{3.8}
\end{align}

If $\delta\leq0$, i.e., $E(0)<0$, we choose $\varepsilon=p$ in {\rm\eqref{3.8}} and $b$ small enough such that $b\leq-2E(0)$. Together with {\rm\eqref{2.7}}, we obtain
\begin{align}
\Phi(t)\geq&\left[(p-2)-\left(p-2+\frac{1}{p}\right)\int_{0}^{t}g(s){\rm d}s\right]\int_{0}^{\ell}x u_{x}^2{\rm d}x+a(p-2)\int_{0}^{t}\int_{0}^{\ell}x u_{s}^2{\rm d}x{\rm d}s\nonumber\\
&+p\left(-2E(0)-b\right)\nonumber\\
\geq&\left[(p-2)-\left(p-2+\frac{1}{p}\right)\int_{0}^{t}g(s){\rm d}s\right]\int_{0}^{\ell}x u_{x}^2{\rm d}x+a(p-2)\int_{0}^{t}\int_{0}^{\ell}x u_{s}^2{\rm d}x{\rm d}s\nonumber\\
\geq&a(p-2)\int_{0}^{t}\int_{0}^{\ell}x u_{s}^2{\rm d}x{\rm d}s\nonumber\\
\geq&0.\label{3.9}
\end{align}

If $0<\delta<1$, i.e., $E(0)<\delta d_{1}$, we choose $\varepsilon=(1-\delta)p+2\delta$ and $b=2(\delta d_{1}-E(0))>0$ in {\rm\eqref{3.8}}. Then we get
\begin{equation*}
\begin{aligned}
\Phi(t)\geq&-2p\delta d_{1}+\left[(p-2)-\left(p-2+\frac{1}{(1-\delta)p+2\delta}\right)\int_{0}^{t}g(s){\rm d}s\right]\int_{0}^{\ell}x u_{x}^2{\rm d}x\\
&+\delta(p-2)(g\circ u_{x})(t)+a(p-2)\int_{0}^{t}\int_{0}^{\ell}xu_{s}^2{\rm d}x{\rm d}s.
\end{aligned}
\end{equation*}
By {\rm\eqref{2.7}}, we have
$$(p-2)-\left(p-2+\frac{1}{(1-\delta)p+2\delta}\right)\int_{0}^{t}g(s){\rm d}s\geq\delta(p-2)\left(1-\int_{0}^{t}g(s){\rm d}s\right)$$
and therefore, by {\rm\eqref{3.1}} and, {\rm\eqref{2.6}} we get

\begin{align}
\Phi(t)\geq&-2p\delta d_{1}+\delta(p-2)\bigg[\left(1-\int_{0}^{t}g(s){\rm d}s\right)\int_{0}^{\ell}x u_{x}^2{\rm d}x+(g\circ u_{x})(t)\bigg]\nonumber\\
&+a(p-2)\int_{0}^{t}\int_{0}^{\ell}xu_{s}^2{\rm d}x{\rm d}s\nonumber\\
\geq&2p\left(\delta d_{1}-\delta d_{1}\right)+a(p-2)\int_{0}^{t}\int_{0}^{\ell}xu_{s}^2{\rm d}x{\rm d}s\nonumber\\
\geq&0.\label{3.10}
\end{align}

Therefore, combining {\rm\eqref{3.5}}, {\rm\eqref{3.9}}, and {\rm\eqref{3.10}}, we arrive at
$$L(t)L''(t)-\frac{p+2}{4}L'(t)^2\geq0, \quad  a.e.\quad t\in [0, T].$$
Let $T_{0}$ be any number which depends only on $p$, $b$, $\int_{0}^{\ell}xu_{0}^2{\rm d}x$ and $\int_{0}^{\ell}xu_{1}^2{\rm d}x$ as
$$T_{0}>\frac{(p-2+4a)\int_{0}^{\ell}xu_{0}^2{\rm d}x+(p-2)\int_{0}^{\ell}xu_{1}^2{\rm d}x}{2(p-2)b},$$
which fulfills the requirement of $$L'(0)=2\int_{0}^{\ell}x u_{0}u_{1}{\rm d}x+2bT_{0}>0.$$
Then using Lemma \ref{le3.3}, we obtain that $L(t)$ goes to $\infty$ as $t$ tends to some $T^*$ satisfying
\begin{equation}\label{3.09032}
 T^*\leq\frac{4L(0)}{(p-2)L'(0)}=\frac{2(1+aT)\int_{0}^{\ell}xu_{0}^2{\rm d}x+2bT_{0}^2}{(p-2)\int_{0}^{\ell}xu_{0}u_{1}{\rm d}x+(p-2)bT_{0}}.
\end{equation}
 Finally, for fixed $T_{0}$, we choose $T$
as
\begin{equation}\label{3.09033}
T>\frac{4\left(\int_{0}^{\ell}xu_{0}^2{\rm d}x+bT_{0}^2\right)}{2(p-2)b T_{0}-(p-2+4a)\int_{0}^{\ell}xu_{0}^2{\rm d}x-(p-2)\int_{0}^{\ell}xu_{1}^2{\rm d}x}.
\end{equation}
Combing \eqref{3.09032} and \eqref{3.09033}, we get $T>T^*$
and this
 contradicts to our assumption, which finishes our proof.

\begin{remark}\label{re7}
We can see that, when $a=0$, Wu \cite{wu2011} established  blow-up
results under some restrictions on $\int_{0}^{\ell}xu_{0}u_{1}{\rm
d}x$, which are no more needed in this paper. In fact, we use the
the potential well theory and the modified convexity method, which
is different from that in Wu \cite{wu2011}.

\end{remark}

\section{Decay of solutions}\label{s4}

In this section we prove our  decay result. For this purpose,
we need the following lemmas.
\begin{lemma}\label{le4.1}{\rm(\cite[Lemma 4.1]{mm10})}
Under the same assumption as in Theorem \ref{th2.6}, one has $I(u(t))>0$ for all $t\in [0, T_{max})$.
\end{lemma}

{\bf Proof of Theorem \ref{th2.4}.} \ We can refer to { \cite[Lemma
4.2]{mm10}}.

Next, we use the following ``modified" functional
\begin{equation}\label{4.5}
F(t):=E(t)+\varepsilon_{1}\Psi(t)+\varepsilon_{2}\chi(t),
\end{equation}
where $\varepsilon_{1}$ and $\varepsilon_{2}$ are positive constants and
\begin{equation}\label{4.6}
\Psi(t)=\xi(t)\int_{0}^{\ell}xu_{t}u{\rm d}x,
\end{equation}
\begin{equation}\label{4.7}
\chi(t)=-\xi(t)\int_{0}^{\ell}xu_{t}\int_{0}^{t}g(t-s)(u(t)-u(s)){\rm d}s{\rm d}x.
\end{equation}

\begin{lemma}\label{le4.2}
For $\varepsilon_{1}$ and $\varepsilon_{2}$ small enough, we have
\begin{equation}\label{4.8}
\alpha_{1}F(t)\leq E(t)\leq \alpha_{2}F(t)
\end{equation}
holds for two positive constants $\alpha_{1}$ and $\alpha_{2}$.
\end{lemma}
{\bf Proof.}\ Straightforward computations lead to
\begin{equation}
\begin{aligned}
F(t)=&E(t)+\varepsilon_{1}\xi(t)\int_{0}^{\ell}xu_{t}u{\rm d}x-\varepsilon_{2}\xi(t)\int_{0}^{\ell}xu_{t}\int_{0}^{t}g(t-s)\left(u(t)-u(s)\right){\rm d}s{\rm d}x\\
\leq &E(t)+\frac{\varepsilon_{1}}{2}\xi(t)\int_{0}^{\ell}xu_{t}^2{\rm d}x+\frac{\varepsilon_{1}}{2}\xi(t)\int_{0}^{\ell}x u^2{\rm d}x\,
+\frac{\varepsilon_{2}}{2}\xi(t)\int_{0}^{\ell}xu_{t}^2{\rm d}x\nonumber\\
&+\frac{\varepsilon_{2}}{2}\xi(t)\int_{0}^{\ell}x\left(\int_{0}^{t}g(t-s)\left(u(t)-u(s)\right){\rm d}s\right)^2{\rm d}x\\
\leq &E(t)+\frac{\varepsilon_{1}}{2}\xi(t)\int_{0}^{\ell}xu_{t}^2{\rm d}x+\frac{\varepsilon_{1}}{2}\xi(t)\int_{0}^{\ell}x u^2{\rm d}x\,
+\frac{\varepsilon_{2}}{2}\xi(t)\int_{0}^{\ell}xu_{t}^2{\rm d}x\nonumber\\
&+\frac{\varepsilon_{2}}{2}\xi(t)\int_{0}^{\ell}x\int_{0}^{t}g(s){\rm d}s\int_{0}^{t}g(t-s)\left(u(t)-u(s)\right)^2{\rm d}s{\rm d}x\\
\leq &E(t)+\frac{(\varepsilon_{1}+\varepsilon_{2})\xi(t)}{2}\int_{0}^{\ell}xu_{t}^2{\rm d}x+\frac{C_{p}\varepsilon_{1}}{2}\xi(t)\int_{0}^{\ell}x u_{x}^2{\rm d}x\nonumber\\
&+\frac{\varepsilon_{2}}{2}(1-l)\xi(t)\int_{0}^{\ell}\int_{0}^{t}xg(t-s)\left(u(t)-u(s)\right)^2{\rm d}s{\rm d}x\\
\leq &E(t)+\frac{(\varepsilon_{1}+\varepsilon_{2})\xi(t)}{2}\int_{0}^{\ell}xu_{t}^2{\rm d}x+\frac{C_{p}\varepsilon_{1}}{2}\xi(t)\int_{0}^{\ell}x u_{x}^2{\rm d}x+\frac{\varepsilon_{2}}{2}(1-l)C_{p}\xi(t)(g\circ u_{x})(t)\\
\leq&\frac{1}{\alpha_{1}}E(t),
\end{aligned}
\end{equation}
and in the same way, we get
\begin{equation*}
\begin{aligned}
F(t)\geq&E(t)-\frac{(\varepsilon_{1}+\varepsilon_{2})\xi(t)}{2}\int_{0}^{\ell}xu_{t}^2{\rm d}x-\frac{C_{p}\varepsilon_{1}}{2}\xi(t)\int_{0}^{\ell}x u_{x}^2{\rm d}x-\frac{\varepsilon_{2}}{2}(1-l)C_{p}\xi(t)(g\circ u_{x})(t)\\
\geq&\left[\frac{1}{2}-\frac{(\varepsilon_{1}+\varepsilon_{2})\xi(t)}{2}\right]\int_{0}^{\ell}xu_{t}^2{\rm d}x+\,
\left(\frac{1}{2}l-\frac{C_{p}\varepsilon_{1}}{2}\xi(t)\right)\int_{0}^{\ell}x u_{x}^2{\rm d}x\\
&+\left[\frac{1}{2}-\frac{C_{p}}{2}\varepsilon_{2}(1-l)\xi(t)\right](g\circ u_{x})(t)-\frac{1}{p}\int_{0}^{\ell}x|u|^p{\rm d}x\\
\geq&\frac{1}{\alpha_{2}}E(t),
\end{aligned}
\end{equation*}
for $\varepsilon_{1}$ and $\varepsilon_{2}$ small enough.

\begin{lemma}\label{le4.30} {\rm(\cite[Lemma 4.5]{mm10})}
Let $v\in L^{\infty}((0, T);H), v_{x}\in L^{\infty}((0, T);H) $ and $g$ be a continuous function on $[0, T]$ and suppose that $0<\tau<1$ and $r>0$. Then there exists a constant $C>0$ such that
\begin{equation*}
\begin{aligned}
&\int_{0}^{t}g(t-s)\|v_{x}(\cdot, t)-v_{x}(\cdot, s)\|_{H}^2{\rm d}s\\
\leq& C\left(\sup_{0<s<T}\|v(\cdot, s)\|_{H}^2 \int_{0}^{t}g^{1-\tau}(s){\rm d}s\right)^{\frac{r-1}{r-1+\tau}}
\left( \int_{0}^{t}g^r(t-s)\|v_{x}(\cdot, t)-v_{x}(\cdot, s)\|_{H}^2{\rm d}s\right)^{\frac{\tau}{r-1+\tau}}.
\end{aligned}
\end{equation*}
\end{lemma}

\begin{lemma}\label{le4.3} {\rm(\cite[Lemma 4.6]{mm10})}
Let $v\in L^{\infty}((0, T);H), v_{x}\in L^{\infty}((0, T);H) $ and $g$ be a continuous function on $[0, T]$ and suppose that $r>0$. Then there exists a constant
$C>0$ such that
\begin{equation*}
\begin{aligned}
&\int_{0}^{t}g(t-s)\|v_{x}(\cdot, t)-v_{x}(\cdot, s)\|_{H}^2{\rm d}s\\
\leq& C\left(t\|v_{x}(\cdot, t)\|_{H}^2+\int_{0}^{t}\|v_{x}(\cdot, s)\|_{H}^2{\rm d}s\right)^{\frac{r-1}{r}}\left( \int_{0}^{t}g^r(t-s)\|v_{x}(\cdot, t)-v_{x}(\cdot, s)\|_{H}^2{\rm d}s\right)^\frac{1}{r}.
\end{aligned}
\end{equation*}
\end{lemma}

\begin{lemma}\label{le4.4}

Assume that $2<p<3$ and that {\rm (G1)}, {\rm (G2)} and \eqref{2.5} hold. Then the functional $\Psi(t)$, defined by \eqref{4.6}, satisfies
\begin{align}
\Psi'(t)\leq &\left(1+\frac{a}{2\beta}+\frac{L}{2\alpha}\right)\xi(t)\int_{0}^{\ell}xu_{t}^2{\rm d}x-\left(\frac{l-a\beta C_{p}-\alpha C_{p}L}{2}\right)\xi(t)\int_{0}^{\ell}xu_{x}^2{\rm d}x\nonumber\\
&+\frac{\xi(t)}{2l}\left(\int_{0}^{t}g^{2-r}(s){\rm d}s\right)(g^r\circ u_{x})(t)+\xi(t)\|u\|_{L_{x}^p}^p,\label{4.9}
\end{align}
for all $\alpha,\beta>0$.
\end{lemma}

{\bf Proof.}\ By using the differential equation in \eqref{1.1}, we easily see that
\begin{align}
\Psi'(t)=&\xi(t)\int_{0}^{\ell}xu_{t}^2{\rm d}x+\xi(t)\int_{0}^{\ell}xuu_{tt}{\rm d}x+\xi'(t)\int_{0}^{\ell}xuu_{t}{\rm d}x\nonumber\\
=&\xi(t)\int_{0}^{\ell}xu_{t}^2{\rm d}x-\xi(t)\int_{0}^{\ell}xu_{x}^2{\rm d}x+\xi(t)\int_{0}^{\ell}x|u|^p{\rm d}x-a\xi(t)\int_{0}^{\ell}xuu_{t}{\rm d}x\nonumber\\
&+\xi(t)\int_{0}^{\ell}xu_{x}\int_{0}^{t}g(t-s)u_{x}
(x, s){\rm d}s{\rm d}x+\xi'(t)\int_{0}^{\ell}xuu_{t}{\rm d}x.\label{4.10}
\end{align}

By Young's inequality, {\rm (G1)}, {\rm (G2)}, Lemma \ref{le2.1} and direct calculations, we arrive at (see \cite{mm10})
\begin{align}
&\xi(t)\int_{0}^{\ell}xu_{x}\int_{0}^{t}g(t-s)u_{x}(x, s){\rm d}s{\rm d}x\nonumber\\
\leq &\frac{\xi(t)}{2}\int_{0}^{\ell}xu_{x}^2{\rm d}x+\frac{\xi(t)}{2}
\int_{0}^{\ell}x\bigg[\int_{0}^{t}g(t-s)\left(|u_{x}(s)-u_{x}(t)|+|u_{x}(t)|\right){\rm d}s\bigg]^2{\rm d}x\nonumber\\
\leq &\frac{\xi(t)}{2}\int_{0}^{\ell}xu_{x}^2{\rm d}x+\frac{\xi(t)}{2}(1+\eta)(1-l)^2\int_{0}^{\ell}xu_{x}^2{\rm d}x\nonumber\\
&+\frac{\xi(t)}{2}\left(1+\frac{1}{\eta}\right)\int_{0}^{t}g^{2-r}(s){\rm d}s\int_{0}^{\ell}\int_{0}^{t}xg^r(t-s)|u_{x}(s)-u_{x}(t)|^2{\rm d}s{\rm d}x\label{4.11}
\end{align}
for any $\eta>0$.
We also have
\begin{equation}\label{4.12}
\xi'(t)\int_{0}^{\ell}xuu_{t}{\rm d}x\leq \frac{\xi(t)}{2}\left|\frac{\xi'(t)}{\xi(t)}\right|\left(C_{p}\alpha\int_{0}^{\ell}xu_{x}^2{\rm d}x+\frac{1}{\alpha}\int_{0}^{\ell}xu_{t}^2{\rm d}x\right),\ \forall\ \alpha>0,
\end{equation}
and
\begin{equation}\label{4.13}
-a\xi(t)\int_{0}^{\ell}xuu_{t}{\rm d}x\leq\frac{a\beta C_{p}}{2}\xi(t)\int_{0}^{\ell}xu_{x}^2{\rm d}x+\frac{a}{2\beta}\xi(t)\int_{0}^{\ell}xu_{t}^2{\rm d}x.
\end{equation}
Combining \eqref{4.10}-\eqref{4.13}, we arrive at
\begin{equation*}
\begin{aligned}
\Psi'(t)\leq &\left(1+\frac{L}{2\alpha}+\frac{a}{2\beta}\right)\xi(t)\int_{0}^{\ell}xu_{t}^2{\rm d}x-\frac{\xi(t)}{2}\left[1-(1+\eta)(1-l)^2-aC_{p}\beta-\alpha C_{p}L\right]\int_{0}^{\ell}xu_{x}^2{\rm d}x\\
&+\frac{\xi(t)}{2}\left(1+\frac{1}{\eta}\right)\left(\int_{0}^{t}g^{2-r}(s){\rm d}s\right)(g^r\circ u_{x})(t)+\xi(t)\|u\|_{L_{x}^p}^p.
\end{aligned}
\end{equation*}
By choosing $\eta=\frac{l}{1-l}$, $\eqref{4.9}$ is established.

\begin{lemma}\label{le4.5}
Assume  $2<p<3$ and that {\rm (G1)}, {\rm (G2)} and \eqref{2.5} hold. Then  the functional $\chi(t)$, defined by \eqref{4.7}, satisfies
\begin{align}
\chi'(t)\leq &\xi(t)\theta\left[1+C^*+2(1-l)^2\right]\int_{0}^{\ell}xu_{x}^2{\rm d}x\nonumber\\
&+\xi(t)\left[\theta-\int_{0}^{t}g(s){\rm d}s+a\theta+\theta L\right]\int_{0}^{\ell}xu_{t}^2{\rm d}x\nonumber\\
&+\left[\frac{1}{2\theta}+2\theta+\frac{C_{p}+(a+L)C_{p}}{4\theta}\right]\xi(t)\left(\int_{0}^{t}g^{2-r}(s){\rm d}s\right)(g^r\circ u_{x})(t)\nonumber\\
&-\frac{C_{p}}{4\theta}\xi(t)g(0)(g'\circ u_{x})(t),\label{4.14}
\end{align}
 for all $\theta>0$.
\end{lemma}

{\bf Proof.}\ Direct calculations give

\begin{align}
\chi'(t)=&\xi(t)\int_{0}^{\ell}x u_{x}(t)\left(\int_{0}^{t}g(t-s)(u_{x}(t)-u_{x}(s)){\rm d}s\right){\rm d}x\nonumber\\
&-\xi(t)\int_{0}^{\ell}x\left(\int_{0}^{t}g(t-s)(u_{x}(t)-u_{x}(s)){\rm d}s\right)\left(\int_{0}^{t}g(t-s)u_{x}(s)\right){\rm d}x\nonumber\\
&-\xi(t)\int_{0}^{\ell}x|u|^{p-2}u\left(\int_{0}^{t}g(t-s)(u(t)-u(s)){\rm d}s\right){\rm d}x\nonumber\\
&-\xi(t)\int_{0}^{\ell}xu_{t}\int_{0}^{t}g'(t-s)(u(t)-u(s)){\rm d}s{\rm d}x\nonumber\\
&+a\xi(t)\int_{0}^{\ell}xu_{t}\int_{0}^{t}g(t-s)(u(t)-u(s)){\rm d}s{\rm d}x\nonumber\\
&-\xi(t)\int_{0}^{\ell}xu^{2}_{t}\int_{0}^{t}g(t-s){\rm d}s{\rm d}x-\xi'(t)\int_{0}^{\ell}xu_{t}\int_{0}^{t}g(t-s)(u(t)-u(s)){\rm d}s{\rm d}x.\label{4.15}
\end{align}

We now estimate the right hand side of \eqref{4.15}. For $\theta>0$, similar as in \cite{mm10}, we have the estimates of the first to the fourth terms.

The first term
\begin{align}
&\xi(t)\int_{0}^{\ell}x u_{x}(t)\left(\int_{0}^{t}g(t-s)(u_{x}(t)-u_{x}(s)){\rm d}s\right){\rm d}x\nonumber\\
\leq & \theta\xi(t)\int_{0}^{\ell}x u_{x}^2{\rm d}x+\frac{1}{4\theta}\xi(t)\left(\int_{0}^{t}g^{2-r}(s){\rm d}s\right)(g^r\circ u_{x})(t).\label{4.16}
\end{align}
The second term
\begin{align}
&\xi(t)\int_{0}^{\ell}x\left(\int_{0}^{t}g(t-s)(u_{x}(t)-u_{x}(s)){\rm d}s\right)\left(\int_{0}^{t}g(t-s)u_{x}(s)\right){\rm d}x\nonumber\\
\leq & 2\theta(1-l)^2\xi(t)\int_{0}^{\ell}x u_{x}^2{\rm d}x+\left(2\theta+\frac{1}{4\theta}\right)\xi(t)\left(\int_{0}^{t}g^{2-r}(s){\rm d}s\right)(g^r\circ u_{x})(t).\label{4.17}
\end{align}
The third term
\begin{align}
&\xi(t)\int_{0}^{\ell}x|u|^{p-2}u\left(\int_{0}^{t}g(t-s)(u(t)-u(s)){\rm d}s)\right){\rm d}x\nonumber\\
\leq & \theta C^*\xi(t)\int_{0}^{\ell}xu_{x}^2{\rm d}x
+\xi(t)\frac{C_{p}}{4\theta}\left(\int_{0}^{t}g^{2-r}(s){\rm d}s\right)(g^r\circ u_{x})(t),\label{4.20}
\end{align}
where $C^*=\frac{ C_{*}}{3-p}\left(\frac{2p}{l(p-2)}E(0)\right)^{p-2}$. The fourth term
\begin{align}
&-\xi(t)\int_{0}^{\ell}xu_{t}\int_{0}^{t}g'(t-s)(u(t)-u(s)){\rm d}s{\rm d}x\nonumber\\
\leq &\theta\xi(t)\int_{0}^{\ell}xu_{t}^2{\rm d}x-
\frac{g(0)}{4\theta}C_{p}\xi(t)\left(g'\circ u_{x}\right)(t).\label{4.21}
\end{align}
For the fifth term, by Young's inequality and Lemma \ref{2.1}, we have
\begin{align}
&a\xi(t)\int_{0}^{\ell}xu_{t}\int_{0}^{t}g(t-s)(u(t)-u(s)){\rm d}s{\rm d}x\nonumber\\
\leq &a\theta\xi(t)\int_{0}^{\ell}xu_{t}^2{\rm d}x+\frac{aC_{p}}{4\theta}\xi(t)\left(\int_{0}^{t}g^{2-r}(s){\rm d}s\right)\left(g^r\circ u_{x}\right)(t).\label{4.22}
\end{align}
For the seventh term

\begin{align}
&-\xi'(t)\int_{0}^{\ell}xu_{t}\int_{0}^{t}g(t-s)(u(t)-u(s)){\rm d}s{\rm d}x\nonumber\\
\leq &\xi(t)\left|\frac{\xi'(t)}{\xi(t)}\right|\bigg[\theta\int_{0}^{\ell}xu_{t}^2{\rm d}x+\frac{C_{p}}{4\theta}\left(\int_{0}^{t}g^{2-r}(s){\rm d}s\right)(g^r\circ u_{x})(t)\bigg]\nonumber\\
\leq &\theta L\xi(t)\int_{0}^{\ell}xu_{t}^2{\rm d}x+\frac{C_{p}L}{4\theta}\xi(t)\left(\int_{0}^{t}g^{2-r}(s){\rm d}s\right)(g^r\circ u_{x})(t).\label{4.24}
\end{align}

A combination of \eqref{4.15}-\eqref{4.24} yields \eqref{4.14}.

{\bf Proof of Theorem \ref{th2.6}.} \ Since $g$ is continuous and $g(0)>0$, then for any $t_{0}>0$, we have
\begin{equation}\label{4.26}
\int_{0}^{t}g(s){\rm d}s\geq\int_{0}^{t_{0}}g(s){\rm d}s:=g_{0}, \qquad \forall\ t\geq t_{0}.
\end{equation}
By using \eqref{2.4}, \eqref{4.9}, \eqref{4.14} and \eqref{4.26}, we obtain
\begin{align}
F'(t)=&E'(t)+\varepsilon_{1}\Psi'(t)+\varepsilon_{2}\chi'(t)\nonumber\\
=&\frac{1}{2}\left(g'\circ u_{x}\right)(t)-\frac{1}{2}g(t)\int_{0}^{\ell}xu_{x}^2{\rm d}x-a\int_{0}^{\ell}xu_{t}^2{\rm d}x+\varepsilon_{1}\Psi'(t)+\varepsilon_{2}\chi'(t)\nonumber\\
\leq &-\left[a-\varepsilon_{1}\left(1+\frac{a}{2\beta}+\frac{L}{2\alpha}\right)\xi(t)+\varepsilon_{2}\xi(t)\left(g_{0}-\theta(1+L)-a\theta\right)\right]\int_{0}^{\ell}xu_{t}^2{\rm d}x\nonumber\\
&+\varepsilon_{1}\xi(t)\int_{0}^{\ell}x|u|^p{\rm d}x+\left[\frac{1}{2}-\frac{\varepsilon_{2}\xi(0)}{4\theta}C_{p}g(0)\right]\left(g'\circ u_{x}\right)(t)\nonumber\\
&-\bigg\{\frac{\varepsilon_{1}}{2}\left(l-a\beta C_{p}-\alpha C_{p}L\right)-\varepsilon_{2}\theta\left[(1+C^*+2(1-l)^2\right]\bigg\}
\xi(t)\int_{0}^{\ell}xu_{x}^2{\rm d}x\nonumber\\
&+\left\{\frac{\varepsilon_{1}}{2l}+\varepsilon_{2}\left[\frac{1}{2\theta}+2\theta+\frac{C_{p}+(a+L)C_{p}}{4\theta}\right]\right\}
\xi(t)\left(\int_{0}^{t}g^{2-r}(s){\rm d}s\right)(g^r\circ u_{x})(t)\nonumber\\
\leq &-\left[\frac{a}{\xi(0)}-\varepsilon_{1}\left(1+\frac{a}{2\beta}+\frac{L}{2\alpha}\right)+\varepsilon_{2}\left(g_{0}-\theta(1+L)-a\theta\right)\right]\xi(t)\int_{0}^{\ell}xu_{t}^2{\rm d}x\nonumber\\
&+\varepsilon_{1}\xi(t)\int_{0}^{\ell}x|u|^p{\rm d}x+\left[\frac{1}{2}-\frac{\varepsilon_{2}\xi(0)}{4\theta}C_{p}g(0)\right]\left(g'\circ u_{x}\right)(t)\nonumber\\
&-\bigg\{\frac{\varepsilon_{1}}{2}\left(l-a\beta C_{p}-\alpha C_{p}L\right)-\varepsilon_{2}\theta\left[(1+C^*+2(1-l)^2\right]\bigg\}
\xi(t)\int_{0}^{\ell}xu_{x}^2{\rm d}x\nonumber\\
&+\left\{\frac{\varepsilon_{1}}{2l}+\varepsilon_{2}\left[\frac{1}{2\theta}+2\theta+\frac{C_{p}+(a+L)C_{p}}{4\theta}\right]\right\}
\xi(t)\left(\int_{0}^{t}g^{2-r}(s){\rm d}s\right)(g^r\circ u_{x})(t),\label{4.27}
\end{align}
since $0<\xi(t)\leq \xi(0)$ .

When $a>0$, we choose $\alpha$ and $\beta$ so small that $$l-a\beta C_{p}-\alpha C_{p}L>\frac{l}{2}$$ and then choose $\theta$  small enough satisfying
\begin{equation}\label{4.27'}
k_{2}=\frac{\varepsilon_{1}l}{4}-\varepsilon_{2}\theta\left[(1+C^*+2(1-l)^2\right]>0.
\end{equation}
As far as $\alpha$, $\beta$ and $\theta$ are fixed, we then pick $\varepsilon_{1}$ and $\varepsilon_{2}$ so small that \eqref{4.8} and \eqref{4.27'} remain valid and
$$k_{1}=\frac{a}{\xi(0)}-\varepsilon_{1}\left(1+\frac{a}{2\beta}+\frac{L}{2\alpha}\right)+\varepsilon_{2}\left(g_{0}-\theta(1+L)-a\theta\right)>0, $$
$$k_{3}=\frac{1}{2}-\frac{\varepsilon_{2}C_{p}g(0)}{4\theta}\xi(0)-\left\{\frac{\varepsilon_{1}}{2l}+
\varepsilon_{2}\left[\frac{1}{2\theta}+2\theta+\frac{C_{p}+(a+L)C_{p}}{4\theta}\right]\left(\int_{0}^{t}g^{2-r}(s){\rm d}s\right)\right\}>0.$$
Therefore, using the assumption $g'(t)\leq-\xi(t)g^r(t)$ in {\rm (G2)}, we have, for some $\sigma>0$,
\begin{equation}\label{4.28}
F'(t)\leq-\sigma\xi(t)\left[\int_{0}^{\ell}xu_{t}^2{\rm d}x-\int_{0}^{\ell}x|u|^p{\rm d}x+\int_{0}^{\ell}xu_{x}^2{\rm d}x+(g^r\circ u_{x})(t)\right], \quad \forall\ t\geq t_{0}.
\end{equation}

When $a=0$, we choose $\theta$, $\alpha$ so small that $g_{0}-(1+L)\theta>\frac{1}{2}g_{0}$, $l-\alpha C_{p}L>\frac{l}{2}$, and
$$\frac{4\theta\left[1+C^*+2(1-l)^2\right]}{l}<\frac{g_{0}}{2+\frac{L}{\alpha}}.$$
Whence $\theta$
and $\alpha$ are fixed, the choice of $\varepsilon_{1}$ and $\varepsilon_{2}$ satisfying
$$\frac{4\theta\left[1+C^*+2(1-l)^2\right]}{l}\varepsilon_{2}<\varepsilon_{1}<\frac{g_{0}\varepsilon_{2}}{2+\frac{L}{\alpha}}$$
will make
\begin{equation}\label{4.28'}
k_{1}=-\varepsilon_{1}\left(1+\frac{L}{2\alpha}\right)\xi(0)+\varepsilon_{2}\xi(0)\left(g_{0}-\theta(1+L)\right)>0,
\end{equation}
\begin{equation}\label{4.28''}
k_{2}=\frac{\varepsilon_{1}}{2}\left(l-\alpha C_{p}L\right)-\varepsilon_{2}\theta\left[(1+C^*+2(1-l)^2\right]>0.
\end{equation}
We then pick $\varepsilon_{1}$ and $\varepsilon_{2}$ so small that \eqref{4.8}, \eqref{4.28'} and \eqref{4.28''} remain valid and
 $$k_{3}=\frac{1}{2}-\frac{\varepsilon_{2}C_{p}g(0)}{4\theta}\xi(0)- \left\{\frac{\varepsilon_{1}}{2l}+
\varepsilon_{2}\left[\frac{1}{2\theta}+2\theta+\frac{C_{p}+LC_{p}}{4\theta}\right]\left(\int_{0}^{t}g^{2-r}(s){\rm
d}s\right)\right\}>0.$$ We can still get \eqref{4.28}.

Next, as \eqref{4.28} is proved,  we will give the following two
cases according to the different ranges of $r$:

 Case 1.
$r=1$.

By virtue of the choice of $\varepsilon_{1}$, $\varepsilon_{2}$ and $\theta$, we estimate \eqref{4.28} and obtain, for some constant $\alpha>0$,
\begin{equation}\label{4.29}
F'(t)\leq-\alpha\xi(t)E(t), \qquad \forall\ t\geq t_{0}.
\end{equation}
Hence, with the help of the left hand side inequality in \eqref{4.8} and \eqref{4.29}, we find
\begin{equation}\label{4.30}
F'(t)\leq-\alpha\alpha_{1}\xi(t)F(t), \qquad \forall\ t\geq t_{0}.
\end{equation}
A simple integration of \eqref{4.30} over $(t_{0}, t)$ leads to
\begin{equation}\label{4.31}
F(t)\leq F(t_{0})e^{-(\alpha\alpha_{1})\int_{t_{0}}^{t}\xi(s){\rm d}s}, \qquad \forall\ t\geq t_{0}.
\end{equation}
Therefore, \eqref{2.8} is established by virtue of \eqref{4.8} again.

Case 2. $1<r<\frac{3}{2}$.

By using  \eqref{2.30} we get
$$g(t)^{1-r}\geq (r-1)\int_{t_{0}}^{t}\xi(s){\rm d}s+g(t_{0})^{1-r}.$$
For $\forall\ 0<\tau<1$, we further have
$$\int_{0}^{\infty}g^{1-\tau}(s){\rm d}s\leq\int_{0}^{\infty}\frac{1}{\left[\right(r-1)\int_{t_{0}}^{t}\xi(s){\rm d}s+g(t_{0})^{1-r}]^{\frac{1-\tau}{r-1}}}{\rm d}s.$$
For $0<\tau<2-r<1$, we have $\frac{1-\tau}{r-1}>1$. And using the
fact that $\int_{0}^{+\infty}\xi(s){\rm d}s=+\infty$, we obtain
$$\int_{0}^{\infty}g^{1-\tau}(s){\rm d}s<\infty,\quad \forall\  0<\tau<2-r.$$
So Lemma \ref{le4.30} and \eqref{2.5'} yield
$$(g\circ u_{x})(t)\leq C\left(E(0)\int_{0}^{\infty}g^{1-\tau}(s){\rm d}s\right)^{\frac{r-1}{r-1+\tau}}\left(g^r\circ u_{x}\right)^{\frac{\tau}{r-1+\tau}}\leq C\left(g^r\circ u_{x}\right)^{\frac{\tau}{r-1+\tau}}$$
for some positive constant $C$.
Therefore, for any $r_{1}>1$, we arrive at

\begin{align}
E^{r_{1}}(t)\leq&CE^{r_{1}-1}(0)\left(\int_{0}^{\ell}xu_{t}^2{\rm d}x-\int_{0}^{\ell}x|u|^p{\rm d}x+\int_{0}^{\ell}xu_{x}^2{\rm d}x\right)+C\left(g\circ u_{x}\right)^{r_{1}}\nonumber\\
\leq&CE^{r_{1}-1}(0)\left(\int_{0}^{\ell}xu_{t}^2{\rm d}x-\int_{0}^{\ell}x|u|^p{\rm d}x+\int_{0}^{\ell}xu_{x}^2{\rm d}x\right)+C\left(g^r\circ u_{x}\right)^{\frac{\tau r_{1}}{r-1+\tau}}.\label{4.01}
\end{align}
By choosing $\tau=\frac{1}{2}$ and $r_{1}=2r-1$ (hence $\frac{\tau r_{1}}{r-1+\tau}=1$), estimate \eqref{4.01} gives, for some $\Gamma>0$,
\begin{equation}\label{4.02}
E^{r_{1}}(t)\leq \Gamma\left[\int_{0}^{\ell}xu_{t}^2{\rm d}x-\int_{0}^{\ell}x|u|^p{\rm d}x+\int_{0}^{\ell}xu_{x}^2{\rm d}x+\left(g^r\circ u_{x}\right)(t)\right]
\end{equation}
By combining \eqref{4.8}, \eqref{4.28} and \eqref{4.02}, we obtain
\begin{equation}\label{4.03}
F'(t)\leq-\frac{\sigma}{\Gamma}\xi(t)E^{r_{1}}(t)\leq-\frac{\sigma}{\Gamma}\alpha_{1}^{r_{1}}F^{r_{1}}(t)\xi(t), \quad \forall t\geq t_{0}.
\end{equation}
A simple integration of \eqref{4.03} leads to
\begin{equation}\label{4.04}
F(t)\leq C_{1}\left(1+\int_{t_{0}}^{t}\xi(s){\rm d}s\right)^{-\frac{1}{r_{1}-1}},\quad \forall t\geq t_{0}.
\end{equation}
Therefore,
$$\int_{t_{0}}^{\infty}F(t){\rm d}t\leq C_{1}\int_{t_{0}}^{\infty}\frac{1}{\left(1+\int_{t_{0}}^{t}\xi(s){\rm d}s\right)^\frac{1}{r_{1}-1}}{\rm d}t.$$
Since $\frac{1}{r_{1}-1}>1$ and $1+\int_{t_{0}}^{t}\xi(s){\rm d}s\rightarrow+\infty$ as $t\rightarrow+\infty$, we get
\begin{equation}\label{4.000}
\int_{t_{0}}^{\infty}F(t){\rm d}t<\infty.
\end{equation}

In addition, by using \eqref{1.1.1.2} we have
$$tF(t)\leq\frac{C_{1}t}{\left(1+\int_{t_{0}}^{t}\xi(s){\rm d}s\right)^\frac{1}{r_{1}-1}}\leq C_{r}.$$
Therefore, we obtain
\begin{equation}\label{4.004}
\sup_{t\geq t_{0}}tF(t)<+\infty.
\end{equation}
Since $E(t)$ is bounded, we use \eqref{4.8}, \eqref{4.000}and \eqref{4.004} to get
$$\int_{0}^{\infty}F(t){\rm d}t+\sup_{t\geq 0}tF(t)<\infty.$$
Then, by using \eqref{2.5'} and Lemma \ref{le4.3}, we have
\begin{equation*}
\begin{aligned}
(g\circ u_{x})(t)\leq& C_{2}\left(t\|u_{x}(\cdot, t)\|_{H}^2+\int_{0}^{t}\|u_{x}(\cdot, s)\|_{H}^2{\rm d}s\right)^{\frac{r-1}{r}}\left( \int_{0}^{t}g^r(t-s)\|u_{x}(\cdot, t)-u_{x}(\cdot, s)\|_{H}^2{\rm d}s\right)^\frac{1}{r}\\
\leq& C_{2}\left(tF(t)+\int_{0}^{t}F(s){\rm d}s\right)^{\frac{r-1}{r}}\left(g^r\circ u_{x}\right)^\frac{1}{r}\\
\leq& C_{3}\left(g^r\circ u_{x}\right)^\frac{1}{r},
\end{aligned}
\end{equation*}
which implies that
\begin{equation}\label{4.32}
\left(g^r\circ u_{x}\right)(t)\geq C_{4}(g\circ u_{x})^r,
\end{equation}
for some constant $C_{4}>0$.

Consequently, a  combination of  \eqref{4.28} and \eqref{4.32} yields
\begin{equation*}
F'(t)\leq-C_{5}\xi(t)\left[\int_{0}^{\ell}xu_{t}^2{\rm d}x-\int_{0}^{\ell}x|u|^p{\rm d}x+\int_{0}^{\ell}xu_{x}^2{\rm d}x+(g\circ u_{x})^r\right], \quad \forall\ t\geq t_{0},
\end{equation*}
for some constant $C_{5}>0$.

On the other hand, as in \cite{bm06}, we can get
\begin{equation*}
E^r(t)\leq C_{6}\left[\int_{0}^{\ell}xu_{t}^2{\rm d}x-\int_{0}^{\ell}x|u|^p{\rm d}x+\int_{0}^{\ell}xu_{x}^2{\rm d}x+(g\circ u_{x})^r\right]
\end{equation*}
for all $t\geq0$ and some constant $C_{6}>0$.
Combining the last two inequalities and \eqref{4.8}, we obtain
\begin{equation}\label{4.33}
F'(t)\leq-C_{7}\xi(t)F^r(t), \qquad \forall\ t\geq t_{0}
\end{equation}
for some constant $C_{7}>0$.
A simple integration of \eqref{4.33} over $(t_{0}, t)$ gives
$$F(t)\leq C_{8}\left(1+\int_{t_{0}}^{t}\xi(s){\rm d}s\right)^{-\frac{1}{r-1}}, \qquad \forall\ t\geq t_{0}.$$
Therefore, \eqref{2.8} is obtained by virtue of \eqref{4.8} again.


\subsection*{Acknowledgments}
This work was partly supported by the Tianyuan Fund of Mathematics (Grant No. 11026211), the Natural
Science Foundation of the Jiangsu Higher Education Institutions (Grant No. 09KJB110005) and the JSPS Innovation Program (Grant No. CXLX12\_0490).

\end{document}